\documentclass{amsart} 
\usepackage{amssymb, amsmath}
\usepackage[dvips]{graphicx}
\usepackage{amsfonts}
\usepackage{latexsym}
\usepackage{color}
\usepackage{pifont}
\newtheorem{theorem}{Theorem}	
\newtheorem{lemma}{Lemma}[section]

\newtheorem{fact}{Fact}[section]
\title[Spherical Separation Theorem]
{
Spherical Separation Theorem
}
\author{
Huhe Han}
\address{College of Science, Northwest Agriculture and Forestry University, China}
\email{han-huhe@nwafu.edu.cn}
\author{Takashi Nishimura
}
\date{\today}
\address{Research Group of Mathematical Sciences,  
Research Institute of Environment and Information Sciences,  
Yokohama National University, 
Yokohama 240-8501, Japan}
\email{nishimura-takashi-yx@ynu.ac.jp}
\begin{document}
\begin{abstract}
In this paper, it is shown that for any two non-empty   
closed (resp., open)  and spherical convex subsets 
$\mathcal{W}_1, \mathcal{W}_2$ of $S^n$, 
the intersection $\mathcal{W}_1\cap \mathcal{W}_2$ 
is empty  
if and only if 
the subset 
$\{P\in S^n\; |\; P\cdot Q>0 
\mbox{ for any } 
Q\in \mathcal{W}_1 \mbox{ and } P\cdot R<0 
\mbox{ for any } R\in \mathcal{W}_2\}$ 
is non-empty, open (resp., closed) and spherical convex.    
\end{abstract}
\maketitle  
\section{Introduction}
Throughout this paper, let $n$ and $ S^n$ 
be a positive integer and  the unit sphere of $\mathbb{R}^{n+1}$
respectively.    
\par 
A subset $\mathcal{W}$ of $S^n$ is said to be 
\textit{hemispherical} if there exists a point $P\in S^n$ 
such that $P\cdot Q>0$ for any $Q\in \mathcal{W}$, 
where $P\cdot Q$ stands for the standard scalar product 
for two vectors $P, Q\in \mathbb{R}^{n+1}$.    
A hemispherical subset $\mathcal{W}$ of $S^n$ is said to be \textit{spherical convex} if for any 
$P, Q\in \mathcal{W}$ and any $t\in [0,1]$ 
the unit vector 
\[
\frac{tP+(1-t)Q}{||tP+(1-t)Q||}
\]
is contained in $\mathcal{W}$.     
Spherical convex sets are spherical counterparts of 
convex sets in $\mathbb{R}^n$.     For more detail on 
spherical convex sets, for instance see 
\cite{survey, nishimurasakemi2}.   
\par 
The purpose of this paper is to show the following:  
\begin{theorem}
\label{theorem 1}
Let $\mathcal{W}_1, \mathcal{W}_2\subset S^n$ be 
two non-empty   
closed (resp., open)  and spherical convex subsets.    
Then, the following (1) and (2) are equivalent.   
\begin{enumerate}
\item[(1)] 
$\mathcal{W}_1\cap \mathcal{W}_2=\emptyset$.   
\item[(2)] The subset consisting of points 
$P\in S^n$ such that $P\cdot Q>0$ for any $Q\in \mathcal{W}_1$ and $P\cdot R<0$ for any 
$R\in \mathcal{W}_2$ is non-empty open (resp., closed) and spherical convex. 
\end{enumerate}   
\end{theorem}
Our motivation for proving Theorem \ref{theorem 1} is coming from  
fruitful discussions with the reviewer of Mathematical Reviews  
for the authors' paper \cite{hannishimura}.    
As shown in the review \cite{reem}, some arguments of \cite{hannishimura} 
are not clear.    Under trying to make clear arguments, the authors noticed 
that Theorem \ref{theorem 1} gives clear proofs for some assertions 
pointed out in \cite{reem}, although Theorem \ref{theorem 1} is much 
stronger than a spherical version of hyperplane separation theorem which 
is needed for clear proofs.       
The authors believe that Theorem \ref{theorem 1} is interesting and  significant in itself.   
Therefore, in order to emphasize Theorem \ref{theorem 1}, 
clear proofs based on Theorem 
\ref{theorem 1} are not given in this paper, they are provided elsewhere.     
\par \bigskip 
Theorem \ref{theorem 1} is proved in Section 2.   
\section{Proof of Theorem 1}\label{section 2}
It is clear that (2) implies (1).   Thus, we concentrate on showing 
that (1) implies (2).    
Since $\mathcal{W}_1, \mathcal{W}_2$ are spherical convex, there exist 
$P_1, P_2\in S^n$ such that $P_i\cdot Q_i>0$ 
for any $Q_i\in \mathcal{W}_i$ and any $i\in \{1, 2\}$.   
For any $P\in S^n$, set $H(P)=\{Q\in S^n\; |\; P\cdot Q\ge 0\}$.    
For any $i\in \{1, 2\}$, let 
$\alpha_{{}_{P_i}}: S^n\backslash H(-P_i)\to P_i+T_{{}_{P_i}}(S^n)$ be the  
central projection with respect to $P_i$, namely, 
\[
\alpha_{{}_{P_i}}(Q)=
\frac{1}{(P_i\cdot Q)}Q \quad \mbox{ for any } Q 
\mbox{ such that } P_i\cdot Q >0; 
\]
where $T_{{}_{P_i}}(S^n)$ is the tangent vector space of $S^n$ at $P_i$ and 
 $P_i+T_{{}_{P_i}}(S^n)$ is the tangent affine space of $S^n$ at $P_i$ defined by 
$\left\{P_i+x\; \left|\; x\in T_{{}_{P_i}}(S^n)\right.\right\}$.    
Moreover,  for any $i\in \{1, 2\}$, define the mapping 
$Id_{{}_{P_i}}: T_{{}_{P_i}}(S^n)\to P_i+T_{{}_{P_i}}(S^n)$ by 
\[
Id_{{}_{P_i}}(x)=P_i+x.    
\]     
Set 
\[
\widetilde{\mathcal{W}}_i=Id_{{}_{P_i}}^{-1}\circ \alpha_{{}_{P_i}}\left(\mathcal{W}_i\right).   
\]
for any $i\in \{1, 2\}$.     Then, both $\widetilde{\mathcal{W}}_1, 
\widetilde{\mathcal{W}}_2$ are convex 
and compact (resp., open) if $\mathcal{W}_1, \mathcal{W}_2$ are 
closed (resp., open) subsets of $S^n$.    
\par 
Next, for any $i\in \{1, 2\}$ and any $\varepsilon\ge 0$, set 
\[
D^n_{{}_{P_i}}(\varepsilon)
=\left\{x\in T_{{}_{P_i}}S^n\; |\; ||x||\le \varepsilon \right\}  
\]
and denote the Minkowski-sum 
of $\widetilde{\mathcal{W}}_i$ and $D^n_{{}_{P_i}}(\varepsilon)$ 
by $\widetilde{\mathcal{X}}_i(\varepsilon)$.    
Namely, 
\[
\widetilde{\mathcal{X}}_i(\varepsilon)
=\widetilde{\mathcal{W}}_i +D^n_{{}_{P_i}}(\varepsilon) 
= \left\{x+y\; \left|\; x\in 
\widetilde{\mathcal{W}}_i, \;y\in D^n_{{}_{P_i}}(\varepsilon) \right.\right\}.    
\]
It is well-known that the Minkowski-sum operation and the convex-hull 
operation are commutative (for instance, see \cite{schneider}).   
Thus, we have the following:  
\begin{fact}\label{fact 1}
For any $\varepsilon \ge 0$, both 
$\widetilde{\mathcal{X}}_1(\varepsilon)$  
and 
$\widetilde{\mathcal{X}}_2(\varepsilon)$ 
are convex 
and compact (resp., open) if both $\mathcal{W}_1, \mathcal{W}_2$ are 
closed (resp., open) subsets of $S^n$.   
\end{fact}  
Moreover, the following is clear by the construction.   
\begin{fact}\label{fact 2}
For any $\varepsilon >0$ and any $i\in \{1,2\}$, 
$\widetilde{\mathcal{W}}_i$ is contained in 
\mbox{\rm int}$\left(\widetilde{\mathcal{X}}_i(\varepsilon)\right)$.  
\end{fact} 
\noindent 
Here, int$\left(\widetilde{\mathcal{X}}_i(\varepsilon)\right)$ stands for 
the set consisting of interior points of 
$\widetilde{\mathcal{X}}_i(\varepsilon)$.    
For any $\varepsilon \ge 0$ and any $i\in \{1,2\}$, set 
\[
\mathcal{X}_i(\varepsilon)=
\alpha_{{}_{P_i}}^{-1}\circ 
Id_{{}_{P_i}}\left(\widetilde{\mathcal{X}}_i(\varepsilon)\right)
\]
Then, from Fact \ref{fact 1}, the following holds.   
\begin{fact}\label{fact 3}
For any $\varepsilon \ge 0$ any $i\in \{1,2\}$,  the following two hold:  
\begin{enumerate}
\item[(1)] The intersection 
$\mathcal{X}_i(\varepsilon)\cap H(-P_i)$ is empty.    
\item[(2)] The set $\mathcal{X}_i(\varepsilon)$  
is spherical convex 
and compact (resp., open) if $\mathcal{W}_i$ is 
a closed (resp., open) subset of $S^n$.   
\end{enumerate}   
\end{fact}  
From Fact \ref{fact 2}, the following holds.   
\begin{fact}\label{fact 4}
For any $\varepsilon >0$ and any $i\in \{1,2\}$, 
$\mathcal{W}_i$ is contained in 
\mbox{\rm int}$\left(\mathcal{X}_i(\varepsilon)\right)$.  
\end{fact} 
Moreover, we have the following:   
\begin{lemma}\label{lemma 1}
There exists a positive number $\varepsilon_0$ such that 
\[
\mathcal{X}_1(\varepsilon)\cap \mathcal{X}_2(\varepsilon)=\emptyset
\]
for any $\varepsilon$ ($0\le \varepsilon \le \varepsilon_0$).   
\end{lemma}
\noindent 
\underline{Proof of Lemma \ref{lemma 1}}.\quad 
Suppose that there does not exist a positive real number 
$\varepsilon_0$ satisfying Lemma \ref{lemma 1}.   
Then, there exists a positive number sequence 
$\varepsilon_j$ $(j=1, 2, \ldots)$ such that the following (a) and (b) are 
satisfied:  
\begin{enumerate}
\item[(a)] $\lim_{j\to 0}\varepsilon_j=0$, 
\item[(b)] 
$\mathcal{X}_1(\varepsilon_j)\cap 
\mathcal{X}_2(\varepsilon_j)\ne \emptyset$.    
\end{enumerate}  
From (b), for any $j\in \mathbb{N}$, there exists at least one point 
$Q_j\in \mathcal{X}_1(\varepsilon_j)\cap 
\mathcal{X}_2(\varepsilon_j)$.    Since the point sequence 
$\left\{Q_j\right\}_{j=1, 2, \ldots}$ is inside $S^n$ and $S^n$ is compact, 
there must exist a convergent subsequence 
$\left\{Q_{{}_{j_k}}\right\}$.    Set 
\[
Q=\lim_{j_k\to \infty}Q_{{}_{j_k}}.   
\]
From (a), the following holds for any $i\in \{1,2\}$:   
\[
\mathcal{X}_i\left(\lim_{j_k\to \infty}\varepsilon_{{}_{j_k}}\right)
= 
\mathcal{X}_i(0)=\mathcal{W}_i.   
\]
Thus, we have $Q\in \mathcal{W}_1\cap \mathcal{W}_2$, which 
contradicts the assumption $\mathcal{W}_1\cap \mathcal{W}_2=\emptyset$.  
Therefore, Lemma \ref{lemma 1} holds.     
\hfill $\Box$
\par 
\medskip 
From now on, take a positive real number $\varepsilon$ satisfying 
$\varepsilon\le \varepsilon_0$ and fix it, where $\varepsilon_0$ is the positive real number in Lemma \ref{lemma 1}; and we continue to prove 
that (1) implies (2) in Theorem \ref{theorem 1}.     
For each $i\in \{1,2\}$, let $\mathcal{Y}_i(\varepsilon)$ denote the 
convex hull of $\mathcal{X}_i(\varepsilon)$ in $\mathbb{R}^{n+1}$.    
By Fact \ref{fact 3} and Lemma \ref{lemma 1}, we have the following:   
\begin{fact}\label{fact 5}
The following two hold:  
\begin{enumerate}
\item[(1)] The intersection 
$\mathcal{Y}_1(\varepsilon)
\cap \mathcal{Y}_2(\varepsilon)$ is empty.    
\item[(2)] For any $i\in \{1,2\}$,  
the set $\mathcal{Y}_i(\varepsilon)$  
is convex 
and compact (resp., open) if $\mathcal{W}_i$ is 
a closed (resp., open) subset of $S^n$.   
\end{enumerate}   
\end{fact}  
From Fact \ref{fact 5}, it follows that the hyperplane separation theorem 
in $\mathbb{R}^{n+1}$ can be applied for 
$\mathcal{Y}_1(\varepsilon), \mathcal{Y}_2(\varepsilon)$.    
For detail on the hyperplane separation theorem, see for instance 
\cite{matousek}.    Hence, the following holds:    
\begin{fact}\label{fact 6}
There exists a hyperplane $H\subset \mathbb{R}^{n+1}$ such that 
$\mathcal{Y}_i(\varepsilon)\subset \mathcal{U}_i$, where 
$\mathbb{R}^{n+1}-H=\mathcal{U}_1\cup \mathcal{U}_2$ and  
$\mathcal{U}_1, \mathcal{U}_2$ are open connected subsets of 
$\mathbb{R}^{n+1}$ satisfying 
$\mathcal{U}_1\cap \mathcal{U}_2=\emptyset$.      
\end{fact}
Let $\mathcal{H}$ be the set consisting of hyperplanes in  
$\mathbb{R}^{n+1}$.       
Moreover, let $\mathcal{S}(\varepsilon)$ be the subset of 
$\mathcal{H}$ consisting of hyperplanes given in Fact \ref{fact 6}.  
Then, by Fact \ref{fact 6}, $\mathcal{S}(\varepsilon)$ is non-empty.   
Let $h: S^n\times \mathbb{R}\to \mathcal{H}$ be the mapping 
defined by 
\[
h(P, r)=\left\{x\in \mathbb{R}^{n+1}\; |\; P\cdot x=r\right\}.    
\]
It is clear that $h$ is a surjective mapping giving 
the quotient topology to $\mathcal{H}$.   
Moreover, $h$ is a $2$ to $1$ mapping and thus 
$S^n\times \mathbb{R}$ is the double covering space of 
the topological space $\mathcal{H}$.   
Let $|\pi|: S^n\times \mathbb{R}\to \mathbb{R}$ be the mapping 
defined by 
\[
|\pi|(P, r)=|r|.     
\]
Then, we have the following:   
\begin{lemma}\label{lemma 2}
\[
\inf |\pi|\left(h^{-1}\left(\mathcal{S}(\varepsilon)\right)\right)=0,
\]
where $\varepsilon (\le \varepsilon_0)$ is a fixed number, and $\varepsilon_0$ is the positive real number in Lemma \ref{lemma 1}.
\end{lemma}
\noindent 
\underline{Proof of Lemma \ref{lemma 2}}.\quad 
Suppose that 
\[
\delta=\inf |\pi|\left(h^{-1}\left(\mathcal{S}(\varepsilon)\right)\right)>0.   
\]
By using the positive real number $\delta$, 
define the linear contracting mapping 
$c_\delta: \mathbb{R}^{n+1}\to \mathbb{R}^{n+1}$ by 
\[
c_\delta(x)=\delta x.   
\]
For each $i\in \{1,2\}$, let the convex hull of the union 
$\mathcal{X}_i(\varepsilon)\cup 
c_\delta\left(\mathcal{X}_i(\varepsilon)\right)$ 
be denoted by $\mathcal{Y}_i(\varepsilon, \delta)$.     
By the constructions of $\mathcal{Y}_i(\varepsilon)$ and 
$\mathcal{Y}_i(\varepsilon, \delta)$, the following trivially holds.   
\begin{fact}\label{fact 7}
For each $i\in \{1,2\}$, 
\[
\mathcal{W}_i\subset \mathcal{Y}_i(\varepsilon)\subset 
\mathcal{Y}_i(\varepsilon, \delta).  
\]
\end{fact}
Moreover, by Fact \ref{fact 5} and the constructions of $\mathcal{Y}_i(\varepsilon)$ and $\mathcal{Y}_i(\varepsilon, \delta)$, 
we have the following:   
\begin{fact}\label{fact 8}
The following two hold:  
\begin{enumerate}
\item[(1)] The intersection 
$\mathcal{Y}_1(\varepsilon, \delta)
\cap \mathcal{Y}_2(\varepsilon, \delta)$ is empty.    
\item[(2)] For any $i\in \{1,2\}$,  
the set $\mathcal{Y}_i(\varepsilon,\delta)$  
is convex 
and compact (resp., open) in $\mathbb{R}^{n+1}$ 
if $\mathcal{W}_i$ is 
a closed (resp., open) subset of $S^n$.   
\end{enumerate}   
\end{fact}  
Thus, $\mathcal{Y}_1(\varepsilon, \delta)$ and 
$\mathcal{Y}_2(\varepsilon, \delta)$ satisfy the assumption of 
the hyperplane separation theorem.    
Hence, we have the following:   
\begin{fact}\label{fact 9}
There exists a hyperplane $H\subset \mathbb{R}^{n+1}$ such that 
$\mathcal{Y}_i(\varepsilon, \delta)\subset \mathcal{U}_i$, where 
$\mathbb{R}^{n+1}-H=\mathcal{U}_1\cup \mathcal{U}_2$ and  
$\mathcal{U}_1, \mathcal{U}_2$ are open connected subsets of 
$\mathbb{R}^{n+1}$ satisfying 
$\mathcal{U}_1\cap \mathcal{U}_2=\emptyset$.      
\end{fact}
Let $\mathcal{S}(\varepsilon, \delta)$ be the subset of 
$\mathcal{H}$ consisting of hyperplanes satisfying Fact \ref{fact 9}.  
Fact \ref{fact 7} implies the following:   
\begin{fact}\label{fact 10}  
\[
\mathcal{S}(\varepsilon, \delta)\subset \mathcal{S}(\varepsilon).    
\]
\end{fact}
Fact \ref{fact 10} yields the following:   
\begin{fact}\label{fact 11}
\[
\delta=\inf |\pi|\left(h^{-1}\left(\mathcal{S}(\varepsilon)\right)\right) 
\le 
\inf |\pi|\left(h^{-1}\left(\mathcal{S}(\varepsilon, \delta)\right)\right).   
\]
\end{fact}
On the other hand, since $\triangle{ABO}$ is a right triangle, we have the following (see Figure \ref{figure1}):   

\begin{figure}[htbp]
  \includegraphics[clip,width=8.0cm]{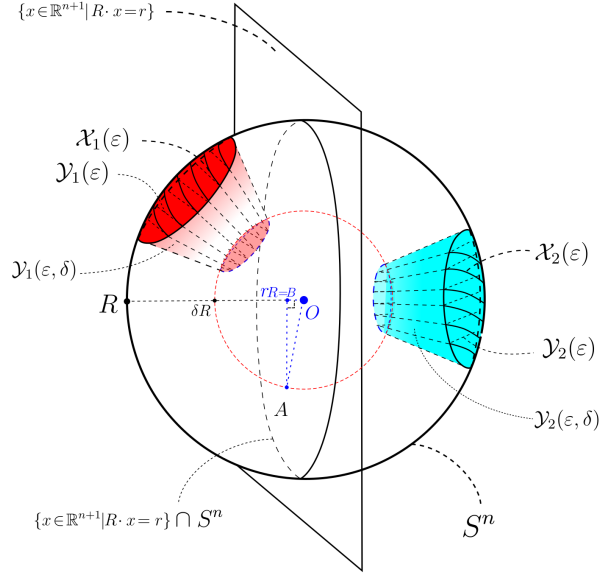}
  \caption{$\mid\pi\mid\left((h^{-1})(\{x\in \mathbb{R}^{n+1}\mid x\cdot R=r\})\right)= |r|<\delta.$}
  \label{figure1}
\end{figure}

\begin{fact}\label{fact 12}
\[
\inf |\pi|\left(h^{-1}\left(\mathcal{S}(\varepsilon, \delta)\right)\right) 
< \delta.   
\]
\end{fact}
Fact \ref{fact 12} contradics Fact \ref{fact 11}.    
Therefore, Lemma \ref{lemma 2} holds.   
\hfill $\Box$
\par 
\medskip 
Fact \ref{fact 4} and Lemma \ref{lemma 2} yield the following:   
\begin{fact}\label{fact 13}
There exists a point $P\in S^n$ such that 
$P\cdot Q > 0$ for any $Q\in \mathcal{W}_1$ and 
$P\cdot R < 0$ for any $R\in \mathcal{W}_2$.     
\end{fact} 
Let $\mathcal{V}$ be the subset of $S^n$ consisting of 
points $P\in S^n$ satisfying  
$P\cdot Q > 0$ for any $Q\in \mathcal{W}_1$ and 
$P\cdot R < 0$ for any $R\in \mathcal{W}_2$.    
Then, Fact \ref{fact 13} implies $\mathcal{V}$ is non-empty.  
\par 
From Fact \ref{fact 5}, Fact \ref{fact 7} and Fact \ref{fact 8}, 
it is easily seen that $\mathcal{V}$ is an open 
(resp., closed) subset of $S^n$ 
if both $\mathcal{W}_1, \mathcal{W}_2$ are closed 
(resp., open) subset of $S^n$.  
The set $\mathcal{V}$ is hemispherical because 
there exists at least one point $Q$ of $\mathcal{W}_1$ 
such that the inequality $P\cdot Q>0$ holds for any 
$P\in \mathcal{V}$, that is to say, 
$\mathcal{V}$ is contained in int$H(Q)$.    
Finally, we show that the set $\mathcal{V}$ is spherical convex.   
This is because 
for any $P_1, P_2\in \mathcal{V}$, 
any $Q\in \mathcal{W}_1$, any  $R\in \mathcal{W}_2$ 
and any $t\in [0,1]$ the following inequalities hold:   
\begin{eqnarray*}
\left(tP_1+(1-t)P_2\right)\cdot Q  =  
t\left(P_1\cdot Q\right)+(1-t)\left(P_2\cdot Q\right) 
& > & 0,  \\ 
\left(tP_1+(1-t)P_2\right)\cdot R  =  
t\left(P_1\cdot R\right)+(1-t)\left(P_2\cdot R\right) 
& < & 0.   
\end{eqnarray*}
\par 
\medskip 
Therefore, (1) implies (2).   
\hfill 
$\Box$


\section*{Acknowledgements}
The authors are grateful to Daniel Reem  
for fruitful discussions. 
The first author was supported by 
the Initial Foundation for Scientific Research of Northwest Agriculture and Forestry University (2452018018).  
The second author was partially supported 
by JSPS KAKENHI Grant Number JP17K05245.   

\end{document}